\documentstyle[12pt,amssymb]{amsart}

\newtheorem{definition}{Definition}
\newtheorem{theorem}{Theorem}
\newtheorem{conjecture}{Conjecture}
\newtheorem{subconjecture}{Sub-Conjecture}

\newcommand{\cB}{{\cal B}}
\newcommand{\C}{{\Bbb C}}

\newcommand{\cE}{{\cal E}}
\newcommand{\cF}{{\cal F}}
\newcommand{\J}{{\cal J}}
\newcommand{\K}{{\cal K}}
\renewcommand{\L}{{\cal L}}

\newcommand{\M}{{\cal M}}
\renewcommand{\O}{{\cal O}}
\renewcommand{\P}{{\Bbb P}}
\newcommand{\Q}{{\Bbb Q}}
\newcommand{\R}{{\Bbb R}}
\renewcommand{\S}{{\Bbb S}}
\newcommand{\Z}{{\Bbb Z}}

\newcommand{\CP}{\C\P}
\newcommand{\dlog}{d\mskip0.5mu\log}
\newcommand{\suchthat}{\ |\ }

\newcommand{\ad}{\operatorname{ad}}
\newcommand{\Aut}{\operatorname{Aut}}
\newcommand{\ch}{\operatorname{ch}}
\newcommand{\dvol}{d\operatorname{vol}}
\newcommand{\Ext}{\operatorname{Ext}}

\newcommand{\Hom}{\operatorname{Hom}}
\renewcommand{\Im}{\operatorname{Im}}
\newcommand{\image}{\operatorname{image}}
\newcommand{\rank}{\operatorname{rank}}
\newcommand{\SU}{\operatorname{SU}}
\newcommand{\supp}{\operatorname{supp}}
\newcommand{\td}{\operatorname{td}}
\newcommand{\U}{\operatorname{U}}
\newcommand{\vol}{\operatorname{vol}}

\begin{document}

\title{The Geometry Underlying Mirror Symmetry}
\author{David R. Morrison}
\address{Department of Mathematics, Box 90320, Duke University,
Durham, NC 27708-0320, USA}
\email{drm@@math.duke.edu}
\begin{abstract}
The recent result of Strominger, Yau and Zaslow relating mirror
symmetry to the quantum field theory notion of T-duality is
reinterpreted as providing a way of geometrically characterizing which
Calabi--Yau manifolds have mirror partners.
The geometric description---that one Calabi--Yau manifold should serve
as a compactified, complexified moduli space for special Lagrangian
tori on the other Calabi--Yau manifold---is rather surprising.  We
formulate some precise mathematical conjectures concerning how these
moduli spaces are to be compactified and complexified, as well as a
definition of geometric mirror pairs (in arbitrary dimension)
which is independent of those
conjectures.  We investigate how this new geometric description ought to
be related to the mathematical statements which have previously been
extracted from mirror symmetry.  In particular, we discuss how the
moduli spaces of the `mirror' Calabi--Yau manifolds should be related
to one another,
and how appropriate subspaces of the homology groups of those manifolds
could be related.  We treat the case of K3 surfaces in some
detail.
\end{abstract}
\maketitle 

Precise mathematical formulations of the string theory phenomenon known as
``mirror symmetry'' \cite{dixon,lvw,CLS,gp} have proved elusive up until
now, largely
due to one of the more mysterious aspects of that symmetry: as
traditionally
formulated, mirror symmetry predicts
an equivalence between physical theories associated to certain pairs
of Calabi--Yau manifolds, but does not specify any geometric relationship
between those manifolds.  However, such a geometric relationship has
recently been discovered in a beautiful paper of Strominger, Yau and Zaslow
\cite{SYZ}.  Briefly put, these authors find that the mirror partner $X$ of
a given Calabi--Yau threefold $Y$
should be realized as the (compactified and
complexified) moduli space for special
Lagrangian tori on $Y$.

This relationship was {derived}\/ in \cite{SYZ} from the assumption
that the physical theories associated to the pair of Calabi--Yau threefolds
satisfy a strong property called ``quantum mirror symmetry''
\cite{S,udual,BBS,mirrorII}.
In the present paper,
we will invert the logic, and use this geometric relationship as a {\it
characterization}\/
 of mirror pairs, which we formulate in arbitrary dimension.\footnote{Our
definition appears to produce valid mirror pairs of conformal field
theories in any dimension, even though the string-theoretic
arguments of \cite{SYZ} cannot be directly extended to
arbitrary dimension
to conclude that all mirror pairs ought to arise in this fashion.}
  On the one hand, this characterization can be stated
in purely mathematical terms, providing a criterion by which mathematicians
can recognize mirror pairs.  On the other hand, the characterization
contains the essential ingredients
needed to apply the quantum field theory
argument known as ``T-duality'' which could in principle
 establish the equivalence of
the associated string theories at the level of physical
rigor (cf.~\cite{towards-duality,SYZ}).\footnote{There are some additional
details which need to be understood before this can be regarded as
fully established in physics.}
This geometric characterization thus
appears to capture the essence of mirror symmetry in mathematical terms.

This paper is organized as follows.  In section 1, we give a brief summary
of quantum mirror symmetry and review the derivation of the geometric
relationship given in \cite{SYZ}.
In section 2, we discuss the theory of special Lagrangian submanifolds
\cite{HL}
and their moduli spaces \cite{McL}, and explain how these moduli spaces
should be compactified and
complexified (following \cite{SYZ}).  In section 3, we review in detail
the topological and Hodge-theoretic properties
which have formed the basis for previous mathematical discussions
of mirror symmetry. We then formulate in section 4
our characterization of {\it geometric mirror pairs}, which we (conjecturally)
relate to those topological and Hodge-theoretic properties.
In
 section 5 we present some new results concerning the geometric
mirror relationship, including a discussion of how it leads to a
connection between certain subspaces
of $H_n(Y)$ and $H_{even}(X)$, and in section 6 we discuss geometric
mirror symmetry for K3 surfaces.

\section{Quantum mirror symmetry}

Moduli spaces which occur in physics  often differ somewhat between the
classical and quantum versions of the same theory.  For example,
the essential mathematical data needed to specify the two-dimensional
conformal field theory
associated to a Calabi--Yau manifold $X$ consists of a Ricci-flat metric
$g_{ij}$ on $X$
and an  $\R/\Z$-valued harmonic $2$-form $B\in{\cal H}^2(X,\R/\Z)$.
The classical version of this theory is independent of $B$ and invariant under
rescaling the metric; one might thus call
the set of all diffeomorphism classes of Ricci-flat metrics of fixed volume
on $X$ the ``classical moduli space'' of the theory.  The volume of the
metric and the $2$-form $B$ must be included in the moduli space
once quantum effects are taken into account; in a
``semiclassical approximation'' to the quantum moduli space, one treats the
data $(g_{ij},B)$ (modulo diffeomorphism) as providing
 a complete description of that  space.  However, a closer analysis of the
physical theory
reveals that this is indeed only an approximation to the quantum moduli space,
with the necessary modifications becoming more and more significant as
the volume is
decreased.  The ultimate source of these modifications---which are of a type
referred to as ``nonperturbative'' in physics---is the set of holomorphic
curves on $X$ and their moduli spaces.
A convenient mathematical way of describing how these modifications work
is this:  there are
certain ``correlation functions'' of the physical theory, which are
described near the large volume limit as power series whose coefficients
are determined by the numbers of holomorphic $2$-spheres on
$X$.\footnote{There are
several possible (equivalent) mathematical interpretations which can
be given to these correlation functions: they can be interpreted
as defining a new ring structure on the cohomology (defining the so-called
quantum cohomology ring) or they can be regarded as defining a variation
of Hodge structure over the moduli space.
We will review
this in more detail in section \ref{s:characterizations} below.}
The quantum moduli
space should then be identified as
the natural domain of definition for these correlation
functions.  To construct it
starting from the semiclassical approximation, one first restricts
to the open set in which the power series converge, and then extends by
analytic continuation to find the complete moduli space.\footnote{There can
also be modifications caused by higher genus curves \cite{chiral}, but these
are less drastic and are not important for our purposes here.}
We refer to this space as the quantum conformal field theory moduli
space $\M_{CFT}(X)$.
(When necessary, we use the notation $\M^{sc}_{CFT}(X)$ to refer to
the semiclassical approximation to this space.)

A similar story has emerged within the last year concerning the moduli
spaces for type IIA and IIB string theories compactified on a Calabi--Yau
threefold $X$.  The classical low-energy physics derived from these string
theories is determined by a quantum conformal field theory, so one might
think of
 the quantum conformal field theory moduli space described above
as being a ``classical moduli space'' for these theories.

In the semiclassical approximation to the {\it quantum}\/ moduli spaces of
these
string theories, we encounter additional mathematical data which must be
specified.
 In the case of the IIA theory, the new data consist of a choice of a
nonzero complex number
(called the ``axion/dilaton expectation value''),
together with an $\R/\Z$-valued harmonic $3$-form $C\in {\cal H}^3(X,\R/\Z)$.
This last object has a familiar mathematical interpretation as a point in
the intermediate Jacobian of $X$ (taking a complex structure on $X$ for which
the metric is K\"ahler).
In the case of the IIB theory on a Calabi--Yau threefold $Y$,
the corresponding new data are a choice of nonzero
``axion-dilaton expectation value''
as before, together with what we might call a quantum
$\R/\Z$-valued harmonic even class $C\in{\cal H}^{\text{even}}_Q(Y,\R/\Z)$.
The word ``quantum'' and the subscript ``$Q$'' here refer to the fact that
we must use the quantum cohomology lattice rather than the ordinary cohomology
lattice in determining when two harmonic $C$'s are equivalent.
(The details
of this difference are not important here; we refer the interested reader
to \cite{udual,mirrorII}.)
For both the IIA and IIB theories, a choice of such ``data'' as above can
be used to describe a low-energy supergravity theory in four dimensions.

Just as in the earlier example, there are additional corrections to the
semiclassical description of the moduli space coming from ``nonperturbative
effects'' \cite{S,GMS,BBS},
among which are some
which go by the name of ``Dirichlet-branes,'' or ``D-branes'' for short.
The source of these D-brane corrections differs for the two string theories
we are considering: in the type IIA theory, they come from
moduli spaces of algebraic cycles on $X$ equipped with flat
$U(1)$-bundles, or more generally, from moduli spaces of
coherent sheaves on $X$.\footnote{A D-brane in type IIA theory
is ordinarily described as a complex submanifold $Z$ together with a flat
$U(1)$-bundle on that submanifold; the associated holomorphic line bundle on
$Z$
can be extended by zero to give a coherent sheaf on $X$.  The {\it arbitrary}\/
coherent sheaves which we consider here
correspond to what are called ``bound states of D-branes'' in physics.
(This same observation has been independently made by Maxim Kontsevich,
and by Jeff Harvey and Greg Moore.)}
In the type IIB theory, the D-brane corrections
come
 from (complexified) moduli spaces of
so-called supersymmetric $3$-cycles on $Y$, the mathematics of which will
be described
in the
next section. Just as the correlation functions which we could use to
determine the structure of the quantum conformal field theory moduli space
 involved a series expansion with contributions from the holomorphic
spheres, the correlation functions in this theory will receive contributions
from the
coherent sheaves
or supersymmetric $3$-cycles, with the precise
nature of the contribution arising from an integral over the
corresponding moduli space.

{\it Quantum mirror symmetry}\/ is the assertion that there should exist
pairs of Calabi--Yau threefolds\footnote{There are versions of quantum
mirror symmetry which can be formulated in other (low) dimensions,
but since these are statements about compactifying ten-dimensional
string theories, they cannot be extended to arbitrarily high dimension.}
$(X,Y)$ such that the type IIA string theory
compactified on $X$ is isomorphic to the type IIB string theory compactified on
$Y$; there should be compatible isomorphisms of both the classical and quantum
theories.  The isomorphism of the classical theories is the statement
that the corresponding (quantum-corrected) conformal field theories should
be isomorphic.  This is the version of mirror symmetry which was translated
into mathematical terms some time ago, and leads to the surprising statements
relating the quantum cohomology on $X$ to the geometric
variation of Hodge structure
on $Y$ (and {\it vice versa}\/).

On the other hand, the isomorphism of the quantum theories has only recently
been explored.\footnote{The speculation some time ago by Donagi and
Markman \cite{DonMark} that some sort of Fourier transform should relate
the continuous data provided by the intermediate Jacobian to the discrete
data provided by the
holomorphic curves is closely related to these isomorphisms
of quantum theories.}
At the semiclassical level, one infers isomorphisms between
the intermediate Jacobian of $X$ (the $3$-form discussed above), and an
analogue of that intermediate Jacobian in quantum cohomology of $Y$.
The full quantum isomorphism would involve properties of the coherent sheaves
on $X$, as related to the supersymmetric $3$-cycles on $Y$.  In fact,
there should be enough correlation functions in the quantum theory to fully
measure the structure of the individual moduli spaces of these sheaves and
cycles,
so we should anticipate that the moduli spaces themselves are
isomorphic.\footnote{In the case of coherent sheaves, one should not use
the usual moduli spaces from algebraic geometry, but rather some sort of
``virtual fundamental cycle'' on the algebro-geometric moduli space, whose
dimension coincides with the ``expected dimension'' of the
algebro-geometric moduli space
as computed from the Riemann--Roch theorem.  When the moduli
problem is unobstructed, this virtual fundamental cycle should coincide with
the usual fundamental cycle on the algebro-geometric moduli space.}
It is this observation which was the key to the Strominger--Yau--Zaslow
argument.

Strominger, Yau and Zaslow observe that the
algebraic $0$-cycles of length one on $X$ (which can be thought of as torsion
sheaves
supported at a single point)
have as their moduli space $X$ itself.  According to quantum mirror
symmetry, then, there should
be a supersymmetric $3$-cycle $M$ on $Y$ with precisely the same moduli space,
that is, the moduli space of $M$ should be $X$.
Since the complex dimension of the moduli space is three, it follows
from a result of McLean \cite{McL} (see the next section) that $b_1(M)=3$.
Now as we will explain in the next section,
the complexified moduli space $\widehat{X}$
for the supersymmetric $3$-cycles parameterizes
both the choice of $3$-cycle $M$ and the choice of a flat $\U(1)$-bundle on
$M$.  Fixing the cycle but varying the bundle gives a real
$3$-torus on $\widehat{X}$ (since $b_1(M)=3$), which turns out to be
a supersymmetric cycle on {\it that}\/ space.  This is the ``inverse''
mirror transform,
based on a cycle $\widetilde{M}$ which is in fact a $3$-torus.
Thus, by applying mirror symmetry
twice if necessary, we see that we can---without loss of
generality---take the original
supersymmetric $3$-cycle $M$ to be a $3$-torus.
In this case, we say that $Y$ has a {\it supersymmetric $T^3$-fibration}; note
that singular fibers must in general be allowed in such fibrations.

We have thus obtained the rough geometric characterization of the pair $(X,Y)$
which was stated in the introduction:
$X$ should be the moduli space for supersymmetric $3$-tori on $Y$.
This characterization
is ``rough'' due to technical difficulties involving both the compactifications
of these moduli spaces, and the complex structures on them.
We will take a different path in section \ref{s:geom}
below, and give a precise geometric characterization which
sidesteps these issues.

This line of argument can be pushed a bit farther, by considering the
 algebraic
$3$-cycle on $X$ in the fundamental class equipped with a flat
$U(1)$-bundle (which must be trivial, and corresponds to the coherent sheaf
$\O_X$).
  There is precisely one of these,
so we find a moduli space consisting of a single point.  Its mirror should
then be a supersymmetric $3$-cycle $M'$ with $b_1(M')=0$.  Moreover, we
should expect quantum mirror symmetry to preserve the intersection theory
of the cycles represented by D-branes
(up to sign),
so since the $0$-cycle and $6$-cycle
on $X$ have intersection number one, we should expect $M$ and $M'$ to have
intersection number one if $M'$ is oriented properly.
In other words, the supersymmetric $T^3$-fibration on $Y$ should have a
section,\footnote{The existence of a section is also expected on
other grounds: the set
of flat $U(1)$-bundles on $M$ has a distinguished element---the trivial
bundle.  This provides a section for the ``dual'' fibration, and suggests
(by a double application of quantum mirror symmetry as above) that the
original fibration could have been chosen to have a section, without
loss of generality.}
 and the base of the fibration should satisfy $b_1=0$.

The final step in the physics discussion given in \cite{SYZ}
is to observe that given a Calabi--Yau threefold with a supersymmetric
$T^3$-fibration and a mirror partner, the mirror partner can be
recovered by dualizing the tori in the fibration, at least generically.
This suggests that by applying an appropriate duality transformation
to the path integral---this
 is known as the ``T-duality argument'' in quantum field theory---one
should be able
to conclude that mirror symmetry does indeed hold for the
corresponding physical theories.
Strominger, Yau and Zaslow
take the first steps towards constructing such an argument, at
appropriate
limit points of the moduli space.  To complete the argument and
extend it to general
points in the moduli space, one would need
to understand the behavior of the T-duality
transformations near the singular fibers; to this end, a detailed
mathematical study of the possible singular fibers is needed.
Some preliminary information about these singularities can be found
in \cite{HL,Br:sub} (see also \cite{GW}).

\section{Moduli of special Lagrangian submanifolds}\label{s:SLAG}

The structure of the supersymmetric $3$-cycles which played a r\^ole
in the previous section was determined in \cite{BBS}, where it was found
that they are familiar mathematical objects known as special Lagrangian
submanifolds.  These are a particular class of submanifolds of Calabi--Yau
manifolds first studied by Harvey and Lawson \cite{HL}.  We proceed to
the definitions.

A {\em Calabi-Yau manifold}\/ is a compact connected orientable
 manifold $Y$ of dimension $2n$
which admits Riemannian metrics whose (global) holonomy is contained in
$\SU(n)$.  For any such metric, there is a complex structure on the
manifold with respect to
which the metric is K\"ahler, and a nowhere-vanishing holomorphic $n$-form
$\Omega$ (unique up to constant multiple).
The complex structure, the $n$-form $\Omega$ and the K\"ahler
form $\omega$ are all covariant constant with respect to the
Levi--Civita connection of the Riemannian metric.  This implies that the metric
is Ricci--flat, and that $\Omega\wedge\overline\Omega$ is a constant
multiple of $\omega^n$.

A {\it special Lagrangian submanifold}\/ of $Y$ is a compact
real $n$-manifold $M$ together
with an immersion
$f:M\to Y$ such that $f^*(\Omega_0)$ coincides with
the induced volume form $\dvol_M$ for
an appropriate choice of holomorphic
$n$-form $\Omega_0$.  Equivalently \cite{HL}, one can require that
(1) $M$ is a Lagrangian submanifold with respect
to the symplectic structure defined by $\omega$, i.e., $f^*(\omega)=0$,
and (2) $f^*(\Im\Omega_0)=0$ for an appropriate $\Omega_0$.
To state this second condition in a way which does not require that $\Omega_0$
be specified, write an arbitrary holomorphic $n$-form $\Omega$ in the
form $\Omega=c\,\Omega_0$,
and note that
\[ \int_Mf^*(\Omega)=c\int_M f^*(\Omega_0)=c\,(\vol M).\]
Thus, the ``appropriate'' $n$-form is given by
\[\Omega_0=\frac{(\vol M)\,\Omega}{\int_Mf^*(\Omega)}\]
and we can replace condition (2) by
\[ \text(2')\qquad f^*\left(\Im\left(\frac{\Omega}{\int_Mf^*(\Omega)}
\right)\right)=0.\]
(The factor of $\vol M$ is a real constant which can be omitted from this
last condition.)

Very few explicit examples of special Lagrangian submanifolds are known.
(This is largely due to our lack of detailed understanding of the Calabi--Yau
metrics themselves.)  One interesting class of examples due to
Bryant \cite{Br:min} comes from Calabi--Yau manifolds which are complex
algebraic varieties defined over the real
numbers: the set of real points on the Calabi--Yau manifold is a special
Lagrangian submanifold.
Another interesting class of examples is the special Lagrangian submanifolds
of a K3 surface, which we will discuss in section 6.

In general, special Lagrangian submanifolds can be deformed, and there will be
a {\it moduli space}\/ which describes the set of all special Lagrangian
submanifolds in a given homology class.  Given a special Lagrangian
$f:M\to Y$ and a deformation of the map $f$, since
$f^*(\omega)=0$, the almost-complex structure on
$Y$ induces a canonical identification between the normal bundle
of $M$ in $Y$
and the tangent bundle of $M$.  Thus, the normal vector field defined by the
deformation can be identified with a $1$-form on $Y$.

The key result concerning the moduli space is due to McLean.

\begin{theorem}[McLean \cite{McL}] \quad

\begin{enumerate}
\item First-order deformations of $f$ are canonically identified
with the space of {\em harmonic}\/ $1$-forms on $Y$.

\item All first-order deformations of $f:M\to Y$ can be extended
to actual deformations.  In particular, the moduli space
$\M _{sL}(M,Y)$ of special Lagrangian maps from $M$ to $Y$ is a smooth
manifold of dimension $b_1(M)$.

\end{enumerate}
\end{theorem}

\noindent
(We have in mind a global structure on $\M_{sL}(M,Y)$ in which two
maps will determine the same point in the moduli space if they
differ by a diffeomorphism of $Y$.)
McLean also observes that $\M=\M _{sL}(M,Y)$ admits a natural
$n$-form $\Theta$ defined by
\[ \Theta(v_1,\dots,v_n)=\int_M\theta_1\wedge\dots\wedge\theta_n\]
where $\theta_j$ is the harmonic $1$-form associated
to $v_j\in T_{\M ,f}$.

As was implicitly discussed in the last section, the moduli spaces of
interest in string theory contain additional pieces of data.  To fully
account for the ``nonperturbative D-brane effects'' in the physical theory
(when $n=3$),
the moduli space which we integrate over must include not only the
choice of special Lagrangian submanifold, but also a choice of flat
$\U(1)$-bundle on that manifold.  If we pick a point $b$ on a manifold $M$,
then the space of flat $\U(1)$-bundles on $M$ is given by
\[\Hom(\pi_1(M,b),\U(1))\cong H^1(M,\R)/H^1(M,\Z).\]
Thus, if we construct a universal family for our special Lagrangian
submanifold problem, i.e., a diagram
\[\begin{array}{cc}
\phantom{\stackrel{f}{\longrightarrow}\quad} {\cal U} \quad
\stackrel{f}{\longrightarrow}
& Y\\[3pt]
\vcenter{\llap{$\scriptstyle p$}}\big\downarrow&\\[4pt]
\M_{sL}(M,Y)&\\
\end{array}\]
with the property that the fibers of $p$ are diffeomorphic to $M$
and $f|_{p^{-1}(m)}$ is the map labeled by $m$, and if $p$ has a section
$s:\M_{sL}(M,Y)\to{\cal U}$, then we can
define a moduli space which includes the data of a flat $\U(1)$-bundle
by
\[\M_D(M,Y):=R^1p_*\R_{\cal U}/R^1p_*\Z_{\cal U},\]
which at each point $m\in\M_{sL}(M,Y)$ specializes to
\[H^1(p^{-1}(m),\R)/H^1(p^{-1}(m),\Z)\cong
\Hom(\pi_1(p^{-1}(m),s(m)),\U(1)).\]
(In the case $n=3$, this is the ``D-brane'' moduli space,
which motivates our notation.)
Note that this space fibers naturally over $\M_{sL}(M,Y)$, and
that there is a section of the fibration, given by the trivial
$\U(1)$-bundles.

Both the base and the fiber of the fibration $\M_D(M,Y)\to\M_{sL}(M,Y)$
have dimension $b_1(M)$, and the fibers are real tori.  In fact,
we expect from the physics that there will be a family of
complex structures on
$\M_D(M,Y)$ making it into a complex manifold of complex dimension
$b_1(M)$.  The real tori should roughly correspond to subspaces obtained
by varying
the arguments of the complex variables while keeping their norms
fixed.  It is expected from the physics that the
complex structure should depend on the choice of both a
Ricci-flat metric on $Y$
and also on an auxiliary harmonic $2$-form $B$.  (This
would make $\M_D(M,Y)$ into a
 ``complexification'' of the moduli space $\M_{sL}(M,Y)$ as
mentioned in the introduction.)  It it not
clear at present precisely how those complex structures are to
be constructed, although in the case $b_1(M)=n$, a method is
sketched in \cite{SYZ} for producing an asymptotic formula
for the Ricci-flat metric which would exhibit the desired dependence
on $g_{ij}$ and $B$, and the first term in that formula is
calculated.\footnote{In the language of \cite{SYZ}, the ``tree-level''
metric on the moduli space is computed, but the instanton corrections
to that tree-level metric are left unspecified.}
The complex structure could in principle be inferred from the metric
if it were known.

Motivated by the Strominger--Yau--Zaslow analysis, we now turn our
attention to the case in which $M$ is an $n$-torus.
The earliest speculations that the special Lagrangian $n$-tori might
play a distinguished r\^ole in studying Calabi--Yau manifolds
were made by McLean \cite{McL}, who pointed
out that if $M=T^n$ then the deformations of $M$ should locally
foliate $Y$.  (There should be no self-intersections nearby since the harmonic
$1$-forms corresponding to the first-order deformations are expected
to have no zeros if the metric on the torus is close to being flat.)
McLean speculated that---by analogy with the K3 case where such elliptic
fibrations are well-understood---if certain degenerations were allowed,
the deformations of $M$ might fill out the whole of $Y$.  We formulate
this as a conjecture (essentially due to McLean).

\begin{conjecture}\label{conj:one}
Suppose that $f:T^n\to Y$ is a special Lagrangian
$n$-torus.  Then there is a natural compactification
$\overline{\M}_{sL}(T^n,Y)$ of the moduli space
$\M_{sL}(T^n,Y)$ and a proper map
$g:Y\to\overline{\M}_{sL}(T^n,Y)$ such that
\[\begin{array}{ccc}
g^{-1}(\M_{sL}(M,Y)) & \hookrightarrow & Y \\[4pt]
\vcenter{\llap{$\scriptstyle g$}}\big\downarrow&&\\[4pt]
\M_{sL}(M,Y)&\\
\end{array}\]
is a universal family of Lagrangian $n$-tori in the same homology class
as $f$.
\end{conjecture}

\begin{definition}
When the properties in conjecture \ref{conj:one} hold, we say that $Y$ has a
{\em special Lagrangian $T^n$-fibration}.
\end{definition}

It is not clear at present what sort of structure should be required of
$\overline{\M}_{sL}(T^n,Y)$: perhaps it should be a manifold with
corners,\footnote{This possibility is suggested by the structure of
toric varieties, the moment maps for which express certain complex
manifolds as $T^n$-fibrations over manifolds with corners (compact convex
polyhedra).}
or perhaps some more exotic singularities should be allowed in the
compactification.  We will certainly want to require that the complex
structures extend to the compactification, and that the section
of the fibration extend to a map $\overline{\M}_{sL}(T^n,Y)\to Y$.

The mirror symmetry analysis of \cite{SYZ} as reviewed in
 the previous section suggests that
the family of dual tori $\M_D(T^n,Y)$ can also be compactified,
resulting in a space which is itself a Calabi--Yau manifold. We formalize
this as a conjecture as well.

\begin{conjecture}\label{conj:two}
The family $\M_D(T^n,Y)$
of dual tori over $\M_{sL}(T^n,Y)$ can be compactified to a manifold
$X$ with a proper map $\gamma:X\to\overline{\M}_{sL}(T^n,Y)$,
such that $X$ admits metrics with $\SU(n)$ holonomy for which the
fibers of $\gamma|_{\gamma^{-1}(\M_{sL}(T^n,Y))}$ are special Lagrangian
$n$-tori.
Moreover, the fibration $\gamma$ admits a section
$\tau:\overline{\M}_{sL}(T^n,Y)\to X$ such that
$\tau(\M_{sL}(T^n,Y))\subset \M_D(T^n,Y)\subset X$ is the zero-section.
\end{conjecture}

\noindent
It seems likely that for an appropriate holomorphic
$n$-form $\Omega_0$ on $X$, the pullback $\tau^*(\Omega_0)$ will coincide
with McLean's $n$-form $\Theta$ when restricted to $\M_{sL}(T^n,Y)$.

The most accessible portion of these conjectures would be the following:

\begin{subconjecture}
The family $\M_D(T^n,Y)$
of dual tori over $\M_{sL}(T^n,Y)$ admits complex structures and
Ricci-flat K\"ahler metrics.  In particular, it has a nowhere vanishing
holomorphic
$n$-form.
\end{subconjecture}

Strominger, Yau and Zaslow have obtained some partial results
concerning this subconjecture, for which we refer the reader to \cite{SYZ}.
It appears, for example, that the construction of the complex structure on
the D-brane moduli space should be local around each torus in the torus
fibration.

\section{Mathematical consequences of mirror symmetry}
\label{s:characterizations}

There is by now quite a long history of extracting mathematical statements
from the physical notion of mirror symmetry.
Many of these work in arbitrary dimension, where there is evidence in
physics for  mirror symmetry among conformal field theories
\cite{gp,higherD}.\footnote{In
low dimension where a string-theory interpretation is possible, this
would become the ``classical'' mirror symmetry which one would also want
to extend to a ``quantum'' mirror symmetry if possible.}
In this section, we review two of those mathematical statements,
presented here as definitions.  As the discussion is a bit technical,
some readers may prefer to skip to the next section, where we
 formulate our new definition
of {\it geometric mirror pairs}\/ inspired by the Strominger--Yau--Zaslow
analysis.  Throughout this section, we let $X$ and $Y$ be Calabi--Yau
manifolds of dimension $n$.

The first prediction one extracts from physics about a mirror pair is
a simple equality of Hodge numbers.

\begin{definition}  We say that the pair
$(X,Y)$ {\em passes  the topological mirror test}\/ if
$h^{n-1,1}(X)=h^{1,1}(Y)$
and $h^{1,1}(X)=h^{n-1,1}(Y)$.
\end{definition}

Many examples of pairs passing this test are known; indeed,
the observation of this ``topological pairing''
 in a class of examples
was one of the initial pieces
of evidence in favor of mirror symmetry \cite{CLS}.
Subsequent
 constructions of Batyrev and Borisov
\cite{batyrev:mirror,borisov,BB:dual}
show that all Calabi--Yau complete intersections in toric
varieties belong to pairs which pass this topological mirror test.

For simply-connected Calabi--Yau threefolds, the Hodge numbers
$h^{1,1}$ and $h^{n-1,1}$ determine all of the others, but in higher
dimension there are more.  Na\"{\i}vely one expects to find that
$h^{p,q}(X)=h^{n-p,q}(Y)$.  However, as was discovered by Batyrev and
collaborators \cite{BatD,BB:Hodge},
the proper interpretation of the numerical invariants of the physical theories
requires a modified notion of ``string-theoretic Hodge numbers''
$h^{p,q}_{st}$;
once this modification has been made, these authors show that
$h^{p,q}_{st}(X)=h^{n-p,q}_{st}(Y)$ for the Batyrev--Borisov pairs
$(X,Y)$ of complete intersections in toric varieties.  The class
of pairs for which this modification is needed
includes some of those given by the Greene--Plesser construction
\cite{gp} for which mirror symmetry of the conformal field theories
has been firmly established in physics, so it would appear that this
modification is truly necessary for a mathematical interpretation
of mirror symmetry.  Hopefully, it too will follow from the geometric
characterization being formulated in this paper.

Going beyond the simple topological properties,
a more precise and detailed prediction arises from identifying the quantum
cohomology of one Calabi--Yau manifold with the geometric variation of Hodge
structure of the mirror partner (in the case that the
Calabi--Yau manifolds have no holomorphic $2$-forms).
We will discuss this prediction in considerable detail, in order to ensure
that this paper has self-contained statements of the conjectures
being proposed within it (particularly the ones in sections \ref{s:geom}
and \ref{s:filtration} below relating the ``old'' and ``new'' mathematical
versions of mirror symmetry).

To formulate this precise prediction,
let $X$ be a Calabi--Yau manifold with $h^{2,0}(X)=0$, and let
$\widetilde\M^{sc}_{CFT}(X)$ be the moduli space of triples $(g_{ij},B,\J)$
modulo diffeomorphism,
where $\J$ is a complex structure for which the metric $g_{ij}$ is
K\"ahler.  The map
$\widetilde\M^{sc}_{CFT}(X)\to\M^{sc}_{CFT}(X)$ is finite-to-one,
so this is another good approximation to the conformal field theory
moduli space.  Moreover, there is a natural map
$\widetilde\M^{sc}_{CFT}(X)\to\M_{cx}(X)$ to the moduli space of
complex structures on $X$, whose fiber over $\J$ is
$\K_{\C}(X_{\J})/\Aut(X_{\J})$,
where
\[\K_{\C}(X_{\J})=\{B+i\,\omega\in {\cal H}^2(X,\C/\Z)\suchthat
\omega\in\K_{\J}\}\]
is the {\it complexified K\"ahler cone\footnote{We are
following the conventions of \cite{icm} rather than those of
\cite{compact,beyond}.} of $X_{\J}$}\/ ($\K_{\J}$ being
its usual K\"ahler cone), and $\Aut(X_{\J})$ is the group of holomorphic
automorphisms of $X_{\J}$.

The moduli space of complex structures $\M_{cx}(X)$ has a variation
of Hodge structure defined on it which is
 of geometric origin: roughly speaking,
one takes a universal family $\pi:{\cal X}\to \M_{cx}(X)$ over the
moduli space and constructs a variation of Hodge structure on
the local system $R^n\pi_*\Z_{\cal X}$ by considering the varying
Hodge decomposition of $H^n(X_{\J},\C)$.  The local system
gives rise to a holomorphic vector bundle
$\cF:=(R^n\pi_*\Z_{\cal X})\otimes \O_{\M_{cx}(X)}$ with a flat connection
$\nabla:\cF\to\Omega^1_{\M_{cx}(X)}\otimes\cF$
(whose flat sections are the sections of the local system), and
the varying Hodge decompositions determine the
{\it Hodge filtration}\/
\[\cF=\cF^0\supset\cF^1\supset\dots\supset\cF^n\supset\{0\},\]
 a filtration
by holomorphic subbundles defined by
\[\cF^p|_{\J}=H^{n,0}(X_{\J})\oplus\dots\oplus H^{p,n-p}(X_{\J}),\]
which is known to satisfy the {\it Griffiths transversality property}\/
\[\nabla(\cF^p)\subset\Omega^1_{\M_{cx}(X)}\otimes\cF^{p-1}.\]
  Conversely,
given the bundle with flat connection and filtration, the complexified local
system
$R^n\pi_*\C_{\cal X}$ can be recovered by taking (local) flat sections,
and the Hodge structures can be reconstructed from the filtration.
The original local system of $\Z$-modules is however additional data,
and cannot be recovered from the bundle, connection and filtration
alone.

The moduli space of complex structures $\M_{cx}(X)$ can be compactified to
a complex space $\overline\M$, to which the bundles $\cF^p$ and the
connection $\nabla$ extend; however, the extended connection $\nabla$
acquires {\it regular singular points}\/ along the boundary $\cB$, which means
that it is a map
\[\nabla:\cF\to\Omega^1(\log\cB)\otimes\cF.\]
The residues of $\nabla$ along boundary components describe the {\it monodromy
transformations}\/ about those components, the same monodromy which defines
the local system.  At normal crossings boundary points there is always an
associated {\it monodromy weight filtration}, which we take to be a filtration
on the homology group $H_n(X)$.

The data of the flat connection and the Hodge filtration
are encoded in the
conformal field theory on $X$ (at least for a sub-Hodge structure
containing $\cF^{n-1}$).\footnote{Note that $\cF^n$ appears directly in the
conformal field theory, and $\cF^{n-1}/\cF^n$ appears as a class
of marginal operators in the conformal field theory.  Thus, the conformal
field theory contains {\it at least}\/ as much of the Hodge-theoretic data
as is described by
the smallest sub-Hodge structure containing $\cF^{n-1}$,
and quite possibly more.}
Since
mirror symmetry reverses the r\^oles of base and fiber in the
map
\[\widetilde\M^{sc}_{CFT}(X)\to\M_{cx}(X),\]
one of the predictions of mirror symmetry will be an isomorphism
between this structure and a similar structure on
$\K_{\C}(X_{\J})/\Aut(X_{\J})$.

In fact, the conformal field theory naturally encodes a variation of
Hodge structure on $\K_{\C}(X_{\J})/\Aut(X_{\J})$.
 To describe this mathematically,
we must choose
a {\it framing}, which is a choice of cone
\[\sigma=\R_+e^1+\dots+\R_+e^r\subset H^2(X,\R)\]
which is generated by a basis $e^1$, \dots, $e^r$ of
$H^2(X,\Z)/\text{torsion}$ and whose interior is contained in
 the K\"ahler cone of $X$.
The complexified K\"ahler part of the semiclassical moduli space
then contains as an open subset the space
\[\M_A(\sigma):=(H^2(X,\R)+i\sigma)/H^2(X,\Z),\]
elements of which can be expanded in the form
$\sum\left(\frac1{2\pi i}\log q_j\right)e^j$, leading to the
alternate description
\[\M_A(\sigma)=\{(q_1,\dots,q_r)\suchthat 0<|q_j|<1\}.\]
The  desired variation
of Hodge structure will be defined on a partial compactification
of this space, namely
\[\overline{\M}_A(\sigma):=\{(q_1,\dots,q_r)\suchthat 0\le|q_j|<1\},\]
which has a distinguished boundary point $\vec{0}=(0,\dots,0)$.

The ingredients we need to define the variation of Hodge structure
are the {\it fundamental Gromov--Witten
invariants\footnote{These can be defined using techniques from
symplectic geometry \cite{Ruan,McDSal,RuanTian} or from algebraic geometry
\cite{KontMan,Kont,BehMan,BehFant,Beh,LiTian}.}
of $X$}, which are trilinear maps
\[\Phi^0_\eta:H^*(X,\Q)\oplus H^*(X,\Q)\oplus H^*(X,\Q)\to\Q.\]
Heuristically, when $A$, $B$ and $C$ are integral classes,
$\Phi^0_\eta(A,B,C)$ should be the number of generically
injective\footnote{We have built the ``multiple cover formula''
\cite{tftrc,manin,voisin} into our definitions.}
holomorphic maps $\psi: \CP^1\to X$ in class $\eta$, such
that $\psi(0)\in Z_A$, $\psi(1)\in Z_B$, $\psi(\infty)\in Z_C$
for appropriate cycles $Z_A$, $Z_B$, $Z_C$ Poincar\'e dual to
the classes $A$, $B$, $C$, respectively.
(The invariants vanish unless $\deg A+\deg B+\deg C=2n$.)
{}From these invariants we can define the {\it Gromov--Witten maps}
$
\Gamma_\eta:H^k(X)\to H^{k+2}(X)
$
by requiring that
\[\Gamma_\eta(A)\cdot B|_{[X]}=\frac{\Phi^0_\eta(A,B,C)}{\eta\cdot C}\]
for $B\in H^{2n-k-2}(X)$, $C\in H^2(X)$.  (This is independent of
the choice of $C$.)

These invariants are usually assembled into the ``quantum cohomology
ring'' of $X$, but here we present this structure in the equivalent
form of a variation of Hodge structure over  $\overline{\M}_A(\sigma)$
degenerating along the boundary.
To do so, we define a holomorphic vector bundle
$\cE:=\left(\bigoplus H^{\ell,\ell}(X)\right)\otimes
\O_{\overline{\M}_A(\sigma)},$
and a flat\footnote{The
flatness of this connection is equivalent to
 the associativity of the product in quantum cohomology.}
 connection $\nabla:\cE\to\Omega^1_{\overline\M}(\log\cB)\otimes\cE$
with regular singular points along the boundary
$\cB= \overline{\M}_A(\sigma)-\M_A(\sigma)$
by the formula\footnote{I am indebted to P. Deligne for advice
\cite{deligne} which led to this form of the formula.}
\[\nabla:=\frac1{2\pi i}\,\left(
\sum_{j=1}^r  \dlog q_j\otimes\ad(e^j)+
\sum_{0\ne\eta\in H_2(X,\Z)}
\dlog\left(\frac1{1-q^{\eta}}\right)\otimes\Gamma_\eta
\right)\]
where $q^\eta=\prod q_j^{e^j(\eta)}$,
and where $\ad(e^j):H^k(X)\to H^{k+2}(X)$ is the adjoint map of the
cup product pairing, defined by $\ad(e^j)(A) = e^j\cup A$.
We also define a ``Hodge filtration''
\[\cE^p:=\left(\bigoplus_{0\le\ell\le n-p}
H^{\ell,\ell}(X)\right)\otimes\O_{\overline{\M}_A(\sigma)},\]
which satisfies $\nabla(\cE^p)\subset \Omega^1_{\overline{\M}}(\log\cB)
\otimes\cE^{p-1}.$
This describes a structure we call the {\em framed $A$-variation of
Hodge structure with framing $\sigma$}.  To be a bit more precise,
we should refer to this as a ``formally degenerating variation of Hodge
structure,''
since the series used to define $\nabla$ is only formal.
(More details about such structures
can be found in \cite{parkcity}; cf.~also \cite{deligne}.) There are also some
subtleties
about passing from a local system of complex vector spaces to a local
system of $\Z$-modules which we shall discuss in section \ref{s:filtration}
below.

The residues of $\nabla$ along the boundary components $q_j=0$ are
the adjoint maps $\ad(e^j)$; the corresponding monodromy weight filtration
at $\vec0$ is simply
\[H_{0,0}(X)\subseteq H_{0,0}(X)\oplus H_{1,1}(X)
\subseteq\dots\subseteq(H_{0,0}(X)\oplus\dots\oplus H_{n,n}(X)).\]
Under mirror symmetry, this maps to the geometric monodromy weight filtration
at an appropriate ``large complex structure limit'' point in
$\overline{\M}_{cx}$
(see \cite{predictions} and references therein).
Note that the class of the $0$-cycle is the monodromy-invariant class in
$H_{even}(X)$; thus, its mirror $n$-cycle will be the monodromy-invariant class
in $H_n(Y)$.

Although the choice of a ``framing'' may look unnatural, the relationship
between different choices of framing is completely understood \cite{compact}
(modulo a conjecture about the action of the automorphism group on the
K\"ahler cone).  Varying the framing corresponds to varying which boundary
point in the moduli space one is looking at, possibly after blowing up the
original boundary
of the moduli space  in order to find an appropriate compactification
containing the desired boundary point.

We finally come to the definition which contains our precise
Hodge-theoretic mirror prediction from physics.

\begin{definition} \label{def:HT}
Let $X$ and $Y$ be Calabi--Yau manifolds with $h^{2,0}(X)=h^{2,0}(Y)=0$.
The pair $(X,Y)$  {\em passes the Hodge-theoretic mirror test}\/ if there
exists a partial compactification $\overline{\M}_{cx}(Y)$
of the complex structure moduli space of $Y$, a neighborhood
$U\subset\overline{\M}_{cx}(Y)$ of a
boundary point $P$
of $\overline{\M}_{cx}(Y)$,
 a framing $\sigma$
for $H^2(X)$, and a ``mirror map'' $\mu:U\to\M_A(\sigma)$
mapping $P$ to $\vec0$
such that $\mu^*$ induces an isomorphism between
$\cE^{n-1}$ and $\cF^{n-1}$ which extends to an isomorphism between
sub-variations of Hodge structure of the
$A$-variation of Hodge structure with framing $\sigma$, and
the geometric formally degenerating
variation of Hodge structure at $P$.
\end{definition}

The restriction to a sub-variation of Hodge structure (which occurs only when
the
dimension of the Calabi--Yau manifold is greater than three)
seems to be necessary in order to get an integer structure on the local system
compatible with the complex variation of Hodge structure.
(We will return to this issue in section \ref{s:filtration}.)
It seems likely
that this is related to the need to pass to ``string-theoretic Hodge
numbers,''
which may actually be measuring the Hodge numbers of the appropriate
sub-Hodge structures.

The property described in the Hodge-theoretic mirror test
can be recast in terms of using the limiting variation of
Hodge structure on $Y$ to make predictions about enumerative geometry
of holomorphic rational curves on $X$.  In this sense, there is a great deal
of evidence in particular cases
(see \cite{predictions,higherD}
and
the references therein).  There are also some specific connections
which have been found
between the variations of Hodge structure associated to mirror pairs
of theories \cite{summing}, as
well as a recent theorem \cite{givental} which proves that the expected
enumerative properties hold for an important class of Calabi--Yau manifolds.

Note that if $(X,Y)$ passes the Hodge-theoretic mirror test in both
directions, then
it passes the topological mirror test (essentially by definition, since
the dimensions of the moduli spaces are given by the Hodge numbers $h^{1,1}$
and $h^{n-1,1}$).

\section{Geometric mirror pairs}\label{s:geom}

We now wish to translate the Strominger--Yau--Zaslow analysis into
a definition of {\it geometric mirror pairs}\/ $(X,Y)$, which we
formulate in arbitrary dimension.
(As mentioned earlier, the arguments of \cite{SYZ} cannot be applied
to conclude that all mirror pairs arise in this way, but it
seems reasonable to suppose that a T-duality argument---applied
to conformal field theories only---would continue to hold.)
The most straightforward such definition
 would say that
$X$ is the compactification of the complexified moduli space of special
Lagrangian $n$-tori on $Y$.  However, as indicated by our conjectures
of section \ref{s:SLAG}, at present we do not have adequate
technical control over the compactification to see that it is a Calabi--Yau
manifold.  So we make instead an indirect definition, motivated by the
following observation:  if we had such a compactified moduli space $X$, then
for generic
$x\in X$ there would be a corresponding special Lagrangian $n$-torus
$T_x\subset Y$, and we could
define an {\em incidence correspondence}
\[ Z=\text{closure of}\ \{(x,y)\in X\times Y\suchthat y\in T_x\}. \]
By definition, the projection $Z\to X$ would have special Lagrangian $n$-tori
as generic fibers.
As we saw earlier, from the analysis of \cite{SYZ}
it is  expected  that generic fibers of the other projection
$Z\to Y$ will also be special Lagrangian $n$-tori.
Furthermore, we should expect that as we vary the metrics on $X$ and on $Y$,
the fibrations by special Lagrangian $n$-tori can be deformed along
with the metrics.
(In fact, it is these dependencies on parameters which should lead
to a ``mirror map'' between moduli spaces.)
Thus, we will formulate our definition using a
{\it family}\/ of correspondences depending on $t\in U$ for some
(unspecified) parameter space $U$.

\begin{definition} A pair of Calabi--Yau manifolds $(X,Y)$ is
a {\em geometric mirror pair}\/ if there is a parameter
space $U$ such that for each $t\in U$ there exist
\begin{enumerate}
\item a correspondence $Z_t\subset(X\times Y)$ which is
the closure of a submanifold of dimension $3n$,

\item maps $\tau_t:X\to Z_t$ and $\widetilde\tau_t:Y\to Z_t$ which
serve as  sections for the projection maps $Z_t\to X$ and $Z_t\to Y$,
respectively,

\item a Ricci-flat metric $g_{ij}(t)$ on $X$ with respect to which
generic fibers of the projection map $Z_t\to Y$ are special Lagrangian
$n$-tori, and

\item a Ricci-flat metric $\widetilde g_{ij}(t)$ on $Y$ with respect to which
generic fibers of the projection map $Z_t\to X$ are special Lagrangian
$n$-tori.
\end{enumerate}
Moreover, for generic $z\in Z_t$, the fibers through $z$ of the two projection
maps must be canonically dual as tori (with origins specified by
$\tau_t$ and $\widetilde \tau_t$).
\end{definition}

In a somewhat stronger form of the definition, we might require that
$U$ be sufficiently large so that the images of the natural maps
$U\to\M_{Ric}(X)$ and $U\to\M_{Ric}(Y)$ to the moduli spaces of Ricci-flat
metrics on $X$ and on $Y$ are open subsets of the respective moduli spaces.
It is too much to hope that these maps would be surjective.
The best picture we could hope for, in fact, would be a diagram
of the form
\[ \M_{Ric}(X) \supseteq U_X
\stackrel{\pi_X}{\twoheadleftarrow\joinrel\relbar} U
\stackrel{\pi_Y}{\relbar\joinrel\twoheadrightarrow}
U_Y \subseteq \M_{Ric}(Y) \]
in which $U_X\subseteq\M_{Ric}(X)$ and $U_Y\subseteq\M_{Ric}(Y)$
are open subsets (near certain boundary points in a compactification
and contained within the set of metrics for which
the semiclassical approximation is valid).
The fibers of $\pi_X$ will have dimension $h^{1,1}(X)$, and if the
induced map is the mirror map each fiber of $\pi_X$ must essentially be the
set of $B$-fields on $X$, i.e., it must be a deformation of the
real torus $H^2(X,\R/\Z)$.
This is compatible with the approximate formula\footnote{The ``tree-level''
formula
given in \cite{SYZ} is subject to unspecified instanton corrections.}
 in \cite{SYZ} for a
family of metrics on $Y$, produced by varying the $B$-field on $X$.

We expect that
geometric mirror symmetry will be related to the earlier mathematical
mirror symmetry properties in the following way.

\begin{conjecture}\label{conj:three}
If $(X,Y)$ is a geometric mirror pair, then the parameter space $U$
and the data in the definition of the geometric mirror pair
can be chosen so that
\begin{enumerate}
\item\label{en:one} $(X,Y)$ passes the topological
mirror test\footnote{Part (\ref{en:one})
is a consequence of part (\ref{en:four}) if $h^{2,0}(X)=h^{2,0}(Y)=0$.},

\item $\pi_X:U\to\M_{Ric}(X)$ lifts to a generically finite map
$\widetilde\pi_X:U\to\widetilde U_X\subseteq\M^{sc}_{CFT}(X)$,

\item $\pi_Y:U\to\M_{Ric}(Y)$ lifts to a generically finite map
$\widetilde\pi_Y:U\to\widetilde U_Y\subseteq\M^{sc}_{CFT}(Y)$, and

\item\label{en:four} if $h^{2,0}(X)=h^{2,0}(Y)=0$, then there are
boundary points $P\in\overline\M_{cx}(Y)$, $P'\in\overline\M_{cx}(X)$
and framings $\sigma$ of $H^2(X)$ and $\sigma'$ of $H^2(Y)$ with
partial compactifications
$\overline{\widetilde U}_X\subset\M_A(\sigma)\times\overline\M_{cx}(X)$
and
$\overline{\widetilde U}_Y\subset\overline\M_{cx}(Y)\times\M_A(\sigma')$
such that
the
 composite map $(\widetilde\pi_X)_*(\widetilde\pi_Y)^*$ extends to
a map $\mu^{-1}\times\mu'$ which consists of mirror maps in both
directions (in the sense of definition \ref{def:HT}).
 In particular, $(X,Y)$ passes the Hodge-theoretic mirror
test.
\end{enumerate}
\end{conjecture}

\noindent
Even in the case that $h^{2,0}(X)\ne0$, there is an induced
map $(\widetilde\pi_X)_*(\widetilde\pi_Y)^*$ which should coincide with the
mirror map between the moduli spaces.

\bigskip

If $X$ has several birational models $X^{(j)}$, then all of the semiclassical
moduli
spaces $\M^{sc}_{CFT}(X^{(j)})$ give rise to a common conformal field
theory moduli space (see \cite{AGM}, or for a more mathematical account,
\cite{beyond}).
If we follow a path between the large radius limit points of two of these
models,
and reinterpret that path in the mirror moduli space, we find a path which
leads from one large complex structure limit point of $\M_{cx}(Y)$ to another.
On the other hand, the calculation of \cite{small} shows that the
homology class of the torus\footnote{Recall that this is the
monodromy-invariant
cycle.} in a special Lagrangian $T^n$-fibration does not change when we
move from one of these regions of $\M_{cx}(Y)$ to another.  Thus, the
moduli space of special Lagrangian $T^n$'s themselves must change as we
move from region to region.    It will be interesting to investigate
precisely how this change comes about.

\section{Mirror cohomology and the
weight filtration} \label{s:filtration}

The ``duality'' transformation which links the two members $X$ and $Y$
of a geometric
mirror pair does not induce any obvious relationship between $H^{1,1}(X)$
and $H^{n-1,1}(Y)$, so it may be difficult to imagine how the
topological mirror test can be passed by a geometric mirror pair.
However, at least for a restricted class of topological cycles,
such a relationship {\it can}\/ be found, as part of a more
general relationship between certain subspaces of $H_{even}(X)$ and
$H_n(Y)$.

Fix a special Lagrangian $T^n$-fibration on $Y$ with a special Lagrangian
section, and consider
$n$-cycles $W\subset Y$
with the property that $W$ is the closure of a submanifold $W_0$
whose intersection with
each nonsingular $T^n$ in the fibration
is either empty, or a sub-torus
of dimension $n-k$ (for some fixed integer $k\le n$).
That is, we assume that $W$ can be generically
described as a $T^{n-k}$-bundle over a $k$-manifold, with the $T^{n-k}$'s
linearly embedded in fibers of the given $T^n$-fibration.  We call such
$n$-cycles {\it pure}.

For any pure $n$-cycle $W\subset Y$, there is a
{\it T-dual cycle}\footnote{In physics, when a  T-duality transformation
is applied to a real torus, a D-brane supported
on a sub-torus is mapped to a D-brane supported on the ``dual'' sub-torus
(of complementary dimension); this can be mathematically identified as the
annihilator.  Here, we apply this principle to a family of sub-tori within
a family of tori.}
$W^\vee\subset X(=\overline\M_D(T^n,Y))$ defined as the closure of
an $n$-manifold $W_0^\vee$ satisfying
\[W_0^\vee\cap (T^n)^*=
\begin{cases}
\text{the annihilator of } W\cap T^n \text{ in } (T^n)^*
&\text{if } W\cap T^n\ne\emptyset\\
\emptyset&\text{otherwise}
\end{cases}
\]
for all smooth fibers $(T^n)^*$ in the dual fibration.
Since the annihilator of an $(n-k)$-torus is a $k$-torus, we see that
$W^\vee$ is generically described as a $T^k$-bundle over a $k$-manifold,
and so it defines  a class in $H_{2k}(X)$.
This is our relationship between the space of pure $n$-cycles on $Y$,
and the even homology on $X$.

Taking the T-duality statements from physics very literally, we are led
to the speculation that pure special Lagrangian $n$-cycles have as their
T-duals certain algebraic cycles on $X$; moreover, the moduli spaces
containing corresponding cycles should be isomorphic.\footnote{As the
referee has pointed out, our ``purity'' condition is probably too strong to
be preserved under deformation, but one can hope that all nearby deformations
of a (pure) special Lagrangian $n$-cycle are reflected in deformations of the
corresponding algebraic cycle.}
(Roughly speaking, the $T^k$-fibration on the corresponding algebraic
$k$-cycle should be given by holding the norms of some system of complex
coordinates on the $k$-cycle fixed, while varying their arguments.)
The simplest cases of this statement we have
already encountered in the Strominger--Yau--Zaslow discussion: the special
Lagrangian $n$-cycles which consist of a single fiber (i.e., $k=0$) are
T-dual to the $0$-cycles of length one on $X$, while a special Lagrangian
$n$-cycle which is the zero-section of the fibration (i.e., $k=n$) is T-dual
to the $2n$-cycle in the fundamental class.  This new construction should
extend that correspondence between cycles to a broader class (albeit still
a somewhat narrow one, since pure cycles are quite special).

In fact, the correspondence should be even broader.  If we begin with an
arbitrary irreducible special Lagrangian $n$-cycle $W$ on $Y$ whose
image in $\M_{sL}(T^n,Y)$ has dimension $k$, then $W$ can be generically
described as a bundle of $(n-k)$-manifolds over the image $k$-manifold.
The T-dual of such a cycle should be a coherent sheaf $\cE$ on $X$ whose
support $Z$
is a complex submanifold of dimension $k$ whose image in $\M_{sL}(T^n,Y)$
is that same $k$-manifold.  Thus, to the homology class of $W$ in $H_n(Y)$
we associate the total homology class in $H_{even}(X)$
of the corresponding coherent
sheaf .\footnote{It appears from both the K3 case discussed in the next
section,
and the analysis of \cite{GHM} that the correct total homology class to use
is the Poincar\'e dual of $\ch({\cal E})\sqrt{\td Y}$.}
Note that since the support has complex dimension $k$, this total homology
class lies in $H_0(X)\oplus H_2(X)\oplus\dots\oplus H_{2k}(X)$.

The homology class of the generic fiber of
$W$ within $T^n$ should determine the sub-tori whose T-duals would
sweep out $Z$; when that homology class is $r$ times a primitive
class, the corresponding coherent sheaf should have generic rank $r$ along $Z$.
For example, a multi-section of the special Lagrangian $T^n$ fibration which
meets the fiber $r$ times should correspond to a coherent sheaf whose support
is all of $X$ and whose rank is $r$.

We have thus found a mapping from the subspace $H_n^{sL}(Y)$ of $n$-cycles
with a special Lagrangian representative, to the subspace
$H_{even}^{alg}(Y)$ of homology classes of algebraic cycles (and coherent
sheaves).  If we consider the Leray filtration on special Lagrangian
$n$-cycles on $Y$
\[\S_k:=\{W\in H_n^{sL}(Y)\suchthat \dim(\image W)\le k\},\]
then this will map to
\[H_0^{alg}(X)\oplus H_2^{alg}(X)\oplus\dots\oplus H_{2k}^{alg}(X)\]
(and the pure $n$-cycles on $Y$
will map to homology classes of algebraic
cycles on $Y$).  But this latter filtration on $H_{even}(X)$ is precisely the
monodromy weight filtration of the $A$-variation of Hodge structures on $X$,
which should be mirror to the geometric
 monodromy weight filtration on $Y$!\footnote{This property of the mapping of
D-branes has also been observed by Ooguri, Oz and Yin \cite{OOY}.}
We are thus
led to the following refinement of conjecture \ref{conj:three}.

\begin{conjecture}\label{conj:four}
If $(X,Y)$ is a geometric mirror pair then there exists
a large complex structure limit point $P\in\overline{\M}_{cx}(Y)$
  corresponding to the mirror
partner $X$, and a sub-variation of the geometric variation of
Hodge structure defined on
$H^{sL}_n(Y)^*$ whose monodromy weight filtration at $P$ coincides
with the Leray filtration for the special Lagrangian $T^n$-fibration on $Y$.
Moreover, under the isomorphism of conjecture \ref{conj:three},
this maps to the sub-variation of the $A$-variation of Hodge structure
defined on $H^{alg}_{even}(X)^*$.
\end{conjecture}

The difficulty in putting an integer structure on the $A$-variation
of Hodge structure stems from the fact that $H^{p,p}(X)$ will in
general not be generated by its intersection with $H^{2p}(X,\Z)$.
However, the {\it algebraic}\/ cohomology $H_{even}^{alg}(X)^*$
does not suffer from this problem: its graded pieces are generated
by integer $(p,p)$-classes.  If conjecture \ref{conj:four} holds,
it explains why there is a corresponding sub-variation of the geometric
variation of  Hodge structure on $Y$, also defined over the integers.
We would thus get corresponding local systems over $\Z$ in addition
to the isomorphisms of complex variations of Hodge structure.

\section{Geometric mirror symmetry for K3 surfaces}

The special Lagrangian submanifolds of a K3 surface can be
studied directly, thanks to the following fact due to Harvey
and Lawson \cite{HL}:  given a Ricci-flat metric on a K3 surface $Y$ and
a special Lagrangian submanifold $M$, there exists a complex
structure on $Y$ with respect to which the metric is K\"ahler, such that
$M$ is a complex submanifold of $Y$.  This allows us to immediately
translate the theory of special Lagrangian $T^2$-fibrations on $Y$
to the standard theory of elliptic fibrations.  In this section,
we will discuss geometric
mirror symmetry for K3 surfaces in some detail.  (Some aspects of
this case have also been worked out by Gross and Wilson \cite{GW},
who went on to study geometric mirror symmetry for the Voisin-Borcea
threefolds of the form $(\mbox{K3}\times T^2)/\Z_2$.)

If we fix a cohomology class $\mu\in H^2(Y,\Z)$ which is
primitive (i.e., $\frac1n\,\mu\not\in H^2(Y,\Z)$ for $1<n\in\Z$)
and satisfies $\mu\cdot\mu=0$, then for any Ricci-flat metric we can find a
compatible
complex structure for which $\mu$ has type $(1,1)$ and
$\kappa\cdot\mu>0$ ($\kappa$ being the K\"ahler form).
The class $\mu$ is then represented by a complex
curve, which moves in a one-parameter family, defining the
structure of an elliptic fibration.  Thus, elliptic fibrations of this
sort exist for every Ricci-flat metric on a K3 surface.\footnote{They
even exist---although possibly in degenerate form---for
the ``orbifold'' metrics which occur at certain
limit points of the moduli space: at those points, $\kappa$ is only
required to be semi-positive, but by the index theorem $\kappa^\perp$
cannot contain an isotropic vector such as $\mu$, so it is still possible
to choose a complex structure such that $\kappa\cdot\mu>0$.}

Our conjecture \ref{conj:one} is easy to verify in this case: as
is well-known, the base of the elliptic fibration on a K3 surface
can be completed
to a $2$-sphere, and the resulting map from K3 to $S^2$ is proper.
In fact, the possible singular fibers are known very explicitly in
this case \cite{kodaira}.

To study conjecture \ref{conj:two}, we need to understand the structure
of the ``complexified'' moduli space $\M_D(T^2,Y)$.
Since a flat $U(1)$-bundle on an elliptic curve is equivalent to
a holomorphic line bundle of degree zero, each point in
 $\M_D(T^2,Y)$ has a natural interpretation as such a bundle
on some particular fiber of the elliptic fibration.  Extending that bundle by
zero, we can regard it
 as a sheaf $\L$ on $Y$, with
$\supp({\L})=\image(f)$.
We thus identify $\M_D(T^2,Y)$ as a
 moduli spaces of such sheaves.

Let us briefly recall the facts about the moduli spaces of simple
sheaves on K3 surfaces, as worked out by Mukai
\cite{Muk:sympl,Muk:bundles}.  First, Mukai showed that for any simple
sheaf $\cE$ on $Y$, i.e., one without any non-constant endomorphisms,
the moduli space $M_{simple}$ is smooth at $[\cE]$ of dimension
$\dim\Ext^1(\cE,\cE)=2-\chi(\cE,\cE).$   (The ``$2$'' in the formula arises
from the spaces
$\Hom(\cE,\cE)$ and $\Ext^2(\cE,\cE)$, each of which has dimension one,
due to the constant endomorphisms in the first case, and their Kodaira--Serre
duals in the second case.)

Second,
Mukai introduced an intersection pairing
on $H^{ev}(Y) = H^0(Y) \oplus H^2(Y) \oplus H^4(Y)$ defined by
\[(\alpha,\beta,\gamma)\cdot(\alpha',\beta',\gamma')=(\beta\cdot\beta'
-\alpha\cdot\gamma'-\gamma\cdot\alpha')|_{[Y]},\]
and a slight modification of the usual Chern character $\ch(\cE)$, defined by
\[v(\cE)=\ch(\cE)\sqrt{\td(Y)}=((\rank(\cE),c_1(\cE),\rank(\cE)+
\frac12(c_1(\cE)^2-2c_2(\cE)),\]
so that the Riemann--Roch theorem reads
\[\chi(\cE,\cF)=v(\cE)\cdot v(\cF).\]
In particular, the moduli space $\M_{simple}(v)$
of simple sheaves with $v(\cE)=v$ has dimension
\[\dim\M_{simple}(v)=2-\chi(\cE,\cE)=2-v\cdot v.\]
In the case of moduli spaces
$\M_{simple}(v)$ of
dimension two, Mukai went on to show that whenever the space is compact,
it must be a K3 surface.

The sheaves $\L$ with support on a curve from our elliptic fibration will
have Mukai class $v(\L)=(0,\mu,0)$ for which $v(\L)\cdot v(\L)=\mu\cdot\mu=0$,
so the moduli space has dimension two.  That is, our moduli space
$\M_D(T^2,Y)$ is contained in $\M_{simple}(0,\mu,0)$ as an open subset.
Our second conjecture will follow if we can show that
 this latter space is compact, or at least admits a
natural compactification.
Whether this is true or not could in principle depend on
the  choice of Ricci-flat metric on $Y$. If we restrict
to metrics with the property that $Y$ is algebraic when given the compatible
complex structure for which $\mu$
defines an elliptic fibration (this is a dense set within
the full moduli space), then techniques of algebraic geometry can be
applied to this problem.  General results of Simpson \cite{simpson}
imply that on an algebraic K3 surface, the set of semistable sheaves
with a fixed Mukai vector $v$ forms a projective variety.  This applies
to our situation with $v=(0,\mu,0)$, and provides the desired compactification.
It is to be hoped that compactifications such as this exist even for
non-algebraic K3 surfaces.

The Mukai class $v=(0,\mu,0)$ should now be mapped under mirror symmetry to
the class of a zero-cycle, or the corresponding sheaf $\O_P$; that Mukai
class is $(0,0,1)$.  In fact, the mirror map known in physics \cite{stringK3}
does precisely that: given any primitive isotropic vector $v$ in $H^{ev}(Y)$,
there is a mirror map which takes it to the vector $(0,0,1)$.  Moreover, it is
easy to
calculate how this mirror map affects complex structures, by specifying how it
affects Hodge structures: if we put a Hodge structure on $H^{ev}(Y)$ in which
$H^0$ and $H^4$ have been specified as type $(1,1)$, then the corresponding
Hodge structure at the mirror image point has $v^\perp/v$ as its $H^2$.

This is {\it precisely}\/ the relationship between Hodge structures on
$Y$ and on $\M_{simple}(v)$ which was found by Mukai \cite{Muk:bundles}!
We can thus identify geometric mirror symmetry for K3 surfaces (which
associates the moduli spaces of zero-cycles and special Lagrangian $T^2$'s)
with the mirror symmetry previously found in physics.  It is amusing
to note that in establishing this relationship, Mukai used elliptic
fibrations and bundles on them in a crucial way.

As suggested in the previous section, such a mirror transformation should
act on the totality of special Lagrangian $2$-cycles.  In fact, it is
known that for at least some K3 surfaces, there is a {\it Fourier--Mukai
transform}\/ which associates sheaves on $\M_{simple}(v)$
to sheaves on $Y$
\cite{BBH}. The map between their homology classes is precisely
the mirror map.\footnote{One example of this is given by a special Lagrangian
section of the $T^2$-fibration, which will map to the class $(1,0,1)$
which is the Mukai vector of the fundamental cycle (i.e., of the structure
sheaf) on the mirror.}
  Thus, proving that there exists such a Fourier--Mukai
transform  for arbitrary K3 surfaces (even non-algebraic
ones) would establish a version of conjecture \ref{conj:four} in
this case.

\subsection*{Acknowledgments}
It is a pleasure to thank
Robbert Dijkgraaf,
Ron Donagi,
Brian Greene,
Mark Gross,
Paul Horja,
Sheldon Katz,
Greg Moore,
Ronen Plesser,
Yiannis Vlassopoulos,
Pelham Wilson,
Edward Witten,
and especially
Robert Bryant and
Andy Strominger
for useful discussions;
I also thank Strominger for communicating the results of \cite{SYZ}
prior to publication,
and the referee for useful remarks on the first version.
I am grateful to the Rutgers physics department for hospitality and support
during the early stages of this work,
to the organizers of the European Algebraic Geometry Conference at the
University of Warwick where this work was first presented,
and to the Aspen Center for Physics where the writing was completed.
This research was partially supported by the National Science Foundation
under grant DMS-9401447.


\ifx\undefined\bysame
\newcommand{\bysame}{\leavevmode\hbox to3em{\hrulefill}\,}
\fi

\end{document}